\documentclass[12pt]{amsart}

\usepackage{amsmath,amssymb,amsbsy,amsfonts,amsthm,latexsym,
                     amsopn,amstext,amsxtra,euscript,amscd,mathrsfs}
\usepackage{amsmath,mathtools,amssymb,amsbsy,amsfonts,latexsym,amsopn,amstext,cite,
                                               amsxtra,euscript,amscd,bm}
\usepackage{mathtools}
\usepackage{todonotes}
\usepackage{url}
\usepackage[colorlinks,linkcolor=blue,anchorcolor=blue,citecolor=blue,backref=page]{hyperref}
\usepackage{multirow}
\usepackage{float}
\usepackage{color}

\usepackage[np]{numprint}
\npdecimalsign{\ensuremath{.}}

\usepackage{bibentry}

\usepackage[norefs,nocites]{refcheck}

\renewcommand*{\backref}[1]{}
\renewcommand*{\backrefalt}[4]{%
    \ifcase #1 (Not cited.)%
    \or        (p.\,#2)%
    \else      (pp.\,#2)%
    \fi}

\begin{document}

\newtheorem{theorem}{Theorem}
\newtheorem{lemma}[theorem]{Lemma}
\newtheorem{claim}[theorem]{Claim}
\newtheorem{cor}[theorem]{Corollary}
\newtheorem{prop}[theorem]{Proposition}
\newtheorem{definition}{Definition}
\newtheorem{question}[theorem]{Question}
\newtheorem{remark}[theorem]{Remark}
\newcommand{\hh}{{{\mathrm h}}}

\numberwithin{equation}{section}
\numberwithin{theorem}{section}
\numberwithin{table}{section}
\def\ssssum{\mathop{\sum\!\sum\!\sum\!\sum}}
\def\sssum{\mathop{\sum\!\sum\!\sum}}
\def\ssum{\mathop{\sum \sum}}
\def\dsum{\mathop{\sum \sum}}
\def\iint{\mathop{\int\ldots \int}}

\def\sA{\mathscr A}
\def\sB{\mathscr B}
\def\sC{\mathscr C}
\def\sD{\Delta}
\def\sE{\mathscr E}
\def\sF{\mathscr F}
\def\sG{\mathscr G}
\def\sH{\mathscr H}
\def\sI{\mathscr I}
\def\sJ{\mathscr J}
\def\sK{\mathscr K}
\def\sL{\mathscr L}
\def\sM{\mathscr M}
\def\sN{\mathscr N}
\def\sO{\mathscr O}
\def\sP{\mathscr P}
\def\sQ{\mathscr Q}
\def\sR{\mathscr R}
\def\sS{\mathscr S}
\def\sU{\mathscr U}
\def\sT{\mathscr T}
\def\sV{\mathscr V}
\def\sW{\mathscr W}
\def\sX{\mathscr X}
\def\sY{\mathscr Y}
\def\sZ{\mathscr Z}

\def \fS{\mathfrak S} 

\def\squareforqed{\hbox{\rlap{$\sqcap$}$\sqcup$}}
\def\qed{\ifmmode\squareforqed\else{\unskip\nobreak\hfil
\penalty50\hskip1em\null\nobreak\hfil\squareforqed
\parfillskip=0pt\finalhyphendemerits=0\endgraf}\fi}

\newfont{\teneufm}{eufm10}
\newfont{\seveneufm}{eufm7}
\newfont{\fiveeufm}{eufm5}
%
%
%
%
\def\frak#1{{\fam\eufmfam\relax#1}}

\newcommand{\bflambda}{{\boldsymbol{\lambda}}}
\newcommand{\bfmu}{{\boldsymbol{\mu}}}
\newcommand{\bfxi}{{\boldsymbol{\xi}}}
\newcommand{\bfrho}{{\boldsymbol{\rho}}}

\newcommand{\bfalpha}{{\boldsymbol{\alpha}}}
\newcommand{\bfbeta}{{\boldsymbol{\beta}}}
\newcommand{\bfphi}{{\boldsymbol{\varphi}}}
\newcommand{\bfpsi}{{\boldsymbol{\psi}}}
\newcommand{\bftheta}{{\boldsymbol{\vartheta}}}
\newcommand{\cM}{\ensuremath{\mathcal{M}} }
\newcommand{\cN}{\ensuremath{\mathcal{N}} }
\def\fK{Frak K}
\def\fT{Frak{T}}

\def\fA{{Frak A}}
\def\fB{{Frak B}}
\def\fC{\mathfrak{C}}

\def \balpha{\bm{\alpha}}
\def \bbeta{\bm{\beta}}
\def \bgamma{\bm{\gamma}}
\def \blambda{\bm{\lambda}}
\def \bchi{\bm{\chi}}
\def \bphi{\bm{\varphi}}
\def \bpsi{\bm{\psi}}

\def\eqref#1{(\ref{#1})}

\def\vec#1{\mathbf{#1}}
\def \vx{\vec{x}}
\def \vX{\vec{X}}
\def \vu{\vec{u}}
\def \vU{\vec{U}}
\def\SAMNp{\sS_{r,p}(\bfalpha; \cM, \cN)}
\def\SAMJp{\sS_{r,p}(\bfalpha; \cM, \cJ)}
\def\SAMJp{\sS_{r,p}(\bfalpha; \cI, \cJ)}
\def\SABMNp{\sS_{r,p}(\bfalpha, \bfbeta ; \cM, \cN)}
\def\SAMNq{\sS_{r,q}(\bfalpha; \cM, \cN)}
\def\SABMNp{\sS_{r,p}(\bfalpha, \bfbeta ; \cM, \cN)}
\def\SABIJp{\sS_{r,p}(\bfalpha, \bfbeta ; \cI, \cJ)}
\def\SABMJp{\sS_{r,p}(\bfalpha, \bfbeta ; \cM, \cJ)}


\def\cA{{\mathcal A}}
\def\cB{{\mathcal B}}
\def\cC{{\mathcal C}}
\def\cD{{\mathcal D}}
\def\cE{{\mathcal E}}
\def\cF{{\mathcal F}}
\def\cG{{\mathcal G}}
\def\cH{{\mathcal H}}
\def\cI{{\mathcal I}}
\def\cJ{{\mathcal J}}
\def\cK{{\mathcal K}}
\def\cL{{\mathcal L}}
\def\cM{{\mathcal M}}
\def\cN{{\mathcal N}}
\def\cO{{\mathcal O}}
\def\cP{{\mathcal P}}
\def\cQ{{\mathcal Q}}
\def\cR{{\mathcal R}}
\def\cS{{\mathcal S}}
\def\cT{{\mathcal T}}
\def\cU{{\mathcal U}}
\def\cV{{\mathcal V}}
\def\cW{{\mathcal W}}
\def\cX{{\mathcal X}}
\def\cY{{\mathcal Y}}
\def\cZ{{\mathcal Z}}
\newcommand{\rmod}[1]{\: \mbox{mod} \: #1}

\def\cg{{\mathcal g}}

\def\vr{\mathbf r}

\def\e{{\mathbf{\,e}}}
\def\ep{{\mathbf{\,e}}_p}
\def\eq{{\mathbf{\,e}}_q}
\def\em{{\mathbf{\,e}}_m}

\def\Tr{{\mathrm{Tr}}}
\def\Nm{{\mathrm{Nm}}}

 \def\SS{{\mathbf{S}}}

\def\lcm{{\mathrm{lcm}}}
\def\ord{{\mathrm{ord}}}

\def\SKr{S}

\def\({\left(}
\def\){\right)}
\def\fl#1{\left\lfloor#1\right\rfloor}
\def\rf#1{\left\lceil#1\right\rceil}

\def\mand{\qquad \mbox{and} \qquad}

\newcommand{\commN}[2][]{\todo[#1,color=yellow]{N: #2}}
\newcommand{\commI}[2][]{\todo[#1,color=green!60]{I: #2}}
\newcommand{\commII}[2][]{\todo[#1,color=red!60]{I: #2}}

%
%
%
%




\hyphenation{re-pub-lished}

\mathsurround=1pt

\def\bfdefault{b}
\overfullrule=5pt

\def \F{{\mathbb F}}
\def \K{{\mathbb K}}
\def \N{{\mathbb N}}
\def \Z{{\mathbb Z}}
\def \Q{{\mathbb Q}}
\def \R{{\mathbb R}}
\def \C{{\mathbb C}}
\def\Fp{\F_p}
\def \fp{\Fp^*}

\def\Kmnp{\cK_p(m,n)}
\def\Kmnq{\cK_q(m,n)}
\def\Kmnp{\cK_p(m,n)}
\def\Kxmnq{\cK_q(\bfxi; m,n)}
\def\Kxmnp{\cK_p(\bfxi; m,n)}
\def\Kxnumnp{\cK_{\nu,p}(\bfxi; m,n)}
\def\Kxnumnq{\cK_{\nu,q}(\bfxi; m,n)}

\def\Kmn{\cK_p(m,n)}
\def\psmn{\psi_p(m,n)}

\def \xbar{\overline x}
\def\e{{\mathbf{\,e}}}
\def\ep{{\mathbf{\,e}}_p}

\title[Bilinear sums of Kloosterman sums]{Bounds on bilinear sums of Kloosterman sums}

 \author[N. Bag]{ Nilanjan Bag}
 
\address{Harish-Chandra Research Institute, Chhatnag Rd, Jhunsi, Uttar Pradesh 211019}
\email{nilanjanbag@hri.res.in}

 \author[I. E. Shparlinski] {Igor E. Shparlinski}

\address{Department of Pure Mathematics, University of New South Wales,
Sydney, NSW 2052, Australia}
\email{igor.shparlinski@unsw.edu.au}

\begin{abstract}
We use some elementary arguments to obtain a new bound on bilinear sums with 
weighted Kloosterman sums which complements those recently obtained by
E.~Kowalski,  P.~Michel and  W.~Sawin~(2020).  
\end{abstract}


\maketitle

\section{Introduction}
 Let $p$ be a  prime number and let  $\F_p^\times$ denote the multiplicative group of the finite field $\F_p$ of $p$ elements. We always assume that $\F_p$ is represented by the elements $\{0, 1, \ldots, p-1\}$. 
 
We now define $r$-multidimensional Kloosterman sums
$$
\cK_{r,p}(n)=\frac{1}{p^{(r-1)/2}} \sum_{\substack{x_1, \ldots, x_r \in\mathbb{F}_p^{\times}\\ 
x_1 \cdots x_r \equiv n \bmod  p}} \ep(x_1 + \cdots + x_r),
$$
where $\ep(t)=\exp(2\pi i t/p)$ for all $t\in\mathbb{R}$. These are also called {\it hyper-Kloosterman sums\/}. Such sums can be interpreted as the inverse Mellin transformation of powers of Gauss sum, thus can be used to study the distribution of Gauss sums.  By the classical bound of Deligne~\cite{Deligne} we know that
\begin{equation}\label{Del-bound}
\cK_{r,p}(n)\leq r.
\end{equation}



Here we continue to study cancellations between such sums. 
In particular, we consider the bilinear {\it Type-I  sums\/}
$$
\SAMNp=\sum_{m\in\cM}\sum_{n\in\cN}\alpha_m K_{r,p}(mn)
$$
with two sets $\cM, \cN \subseteq  \mathbb{F}_p^{\times}$ and complex weights
$\bfalpha=\(\alpha_m\)_{m\in\cM}$,  and also  more general  {\it Type-II  sums\/}$$
\SABMNp=\sum_{m\in\cM}\sum_{n\in\cN}\alpha_m\beta_n K_{r,p}(mn)
$$
with complex weights $\bfalpha=\(\alpha_m\)_{m\in\cM}$
and $\bfbeta=\(\beta_n\)_{n\in\cN}$.

The case  when one or both of the sets $\cM$ and $\cN$ are intervals of  $M$ and $N$ consecutive integers 
$$
\cI = \{A+1, \ldots, A+M\} \subseteq \F_p^\times  \quad \text{and}\quad  \cJ = \{B+1, \ldots, B+N\} \subseteq \F_p^\times
$$ 
with some integers $A, B, M, N$ has received most of attention. 

In the classical case of $r=2$,  when $\cN =\cJ$, for Type-I sums the strongest bound  in the range when $M,N = p^{1/2+o(1)}$, which crucial for   applications to moments of various $L$-functions, is given in~\cite{Shp}, see also~\cite{BFKMM1,BFKMM2,FKM,LSZ1,LSZ2,ShpZha,WuXi,Xi-FKM}. 
For example, one of the applications of this bound in an improvement in~\cite[Theorem~3.1]{Shp} of the error term $p^{-1/68+o(1)}$ of Blomer,  Fouvry,  Kowalski,  Michel and  Mili{\'c}evi{\'c}~\cite[Theorem~1.2]{BFKMM1} in the asymptotic formula  for mixed moments of L-series associated with Hecke eigenforms down to $p^{-1/64+o(1)}$.

On the other hand,  if both sets are intervals,  $\cM =\cI$ and  $\cN =\cJ$
the widest range in which there exists a nontrivial bound, for any $r\ge 2$ is due to 
Kowalski,    Michel  and Sawin~\cite{KMS2}, see also~\cite{FKM,KMS1} for previous results. For example,
 the bound 
of~\cite[Theorem~1.2]{KMS2}
improves~\eqref{Del-bound} for $M, N \ge p^{3/8+\varepsilon}$ for any fixed $\varepsilon>0$, while~\cite[Theorem~1.3]{KMS2} does so for $M, N \ge p^{1/3+\varepsilon}$.  In fact, in both~\cite[Theorem~1.2 and~1.3]{KMS2} the set $\cM$ 
can be more general than an interval, however some restrictions of the size 
of the elements of $\cM$ are still necessary. 

We also notice the method of Shkredov~\cite{Shkr}, which is based on additive combinatorics
and thus also has a potential to be extended to general sets. 

Since the case of an arbitrary set $\cM$ without any restrictions is of independent 
interest, and also has some applications to the average values of the divisor function in 
arithmetic progressions~\cite{KeSh},  Banks and Shparlinski~\cite[Theorem~2.4]{BaSh} 
have extended~\cite[Theorem~1.3]{KMS2}  to arbitrary sets $\cM$ and the weights
 $|\alpha_m|\le 1$, $m \in \cM$, and shown that in this case fixed even integer 
 $\ell\ge 1$, we have
\begin{align*}
|\SAMJp|&\le MN\Bigl(N^{-1/2\ell}
+M^{-1/8\ell}N^{-1/\ell}p^{3/8\ell+1/2\ell^2} \\
&\qquad +M^{-1/2\ell}N^{-1/\ell}p^{1/2\ell+1/2\ell^2}
+N^{-3/2\ell}p^{1/2\ell+1/\ell^2}\Bigr)p^{o(1)}.
\end{align*}
 In particular, the 
bound of~\cite[Theorem~2.4]{BaSh}  remains nontrivial in the same range 
$M, N \ge p^{1/3+\varepsilon}$ as~\cite[Theorem~1.3]{KMS2} 
as for an arbitrary set $\cM$.  

\begin{remark}
We note that comparing with~\cite[Theorem~4.3]{KMS2} the above bound requires 
$\ell$ to be even. It seems that this condition is actually missing in the formulation of~\cite[Theorem~4.3]{KMS2}, as the conclusion that $\gamma \ge 1/2$ at the very end 
of~\cite[Section~4]{KMS2}, where it appeals to~\cite[Theorem~4.5]{KMS2} requires 
the inequality $\rf{(\ell-1)/2} \ge \ell/2$, which fails for odd $\ell$. This is, certainly inconsequential for~\cite[Theorem~1.3]{KMS2} (we also note that~\cite[Theorems~1.2 and~4.1]{KMS2} are not effected). 
\end{remark}

Here we obtain analogous results for $\SABMJp$.

%
%
%

\section{General notation}

We  define the norms
 $$
 \|\bfalpha\|_\infty=\max_{m\in \cI}|\alpha_m|  \mand \|\bfalpha\|_\sigma =\( \sum_{m\in \cI} |\alpha_m|^\sigma\)^{1/\sigma},
 $$
where  $\sigma >0$ and similarly for other sequences.

Since the symbol $\overline{x}$ is reserved for the modular inverse of $x \in \F_p^\times$, we use $\overline{z}$
to denote the complex conjugate of $z \in \C$.

Throughout the paper,  as usual $A\ll B$  is equivalent to the inequality $|A|\le cB$
with some constant $c>0$,
 which occasionally, where obvious, may depend
on the integer parameter $\nu \ge 1$, and is absolute otherwise.
The letter $p$ always denotes a prime number.

For a positive $A$ also write $a \sim A$ to denote that $a$ is in the dyadic interval $A \le a \le 2A$.

\section{Main results}

We remark that when the weights   $\balpha$ satsify $|\alpha_m|\leq 1$, $m \in \cM$, then 
by~\eqref{Del-bound}  and the Cauchy inequality we have 
$$
\left|\SABMNp\right| \le  \|\beta\|_2MN^{1/2}, 
$$
where 
$$
 \|\bbeta\|_2 =\( \sum_{n\in \cJ} |\beta_n|^2\)^{1/2}.
$$

\begin{theorem}\label{thm:SABMN}
Consider complex weights $\balpha$ and $\bbeta$ with $|\alpha_m|\leq 1$. Then for any  fixed integer $\ell\geq 2$, 
for set $\cM\subseteq\mathbb{F}_p^{\times}$ of cardinality $M$ and  an interval 
$\cN =   \{B+1, \ldots, B+N\}\subseteq \F_p^\times$ 
of length $p>N\ge p^{3/{2\ell}}$, we have
$$
 \left|\SABMNp\right| \le   \|\beta\|_2MN^{1/2}\Delta p^{o(1)}
$$  
where 
\begin{align*}
 \Delta & =M^{-1/2}  + M^{-1/{4\ell}}N^{-1/{4\ell}}p^{1/{8\ell}}+M^{-5/{16\ell}}N^{-1/{2\ell}}p^{5/{16\ell}+3/{8{\ell}^2}}  \\
 & \qquad \qquad +M^{-1/{2\ell}}N^{-1/{2\ell}}p^{3/{8\ell}+3/{8{\ell}^2}}+M^{-1/{4\ell}}N^{-3/{4\ell}} p^{3/{8\ell}+3/{4{\ell}^2}}  .
\end{align*}
\end{theorem}

\begin{remark}
Our bound in Theorem~\ref{thm:SABMN} is weaker than the bound in~\cite[Theorem~4.1]{KMS2}, but it applies in 
larger generality without any restrictions on $M^{+} = \max\{m:~m \in \cM\}$, which is given in~\cite{KMS2}.
 \end{remark}

 In particular,  taking a sufficiently large values of $\ell$, after simple calculations we derive from Theorem~\ref{thm:SABMN}:
 
 \begin{cor}\label{cor:SABMN-range} Let $\varepsilon > 0$ be fixed.
Consider complex weights $\balpha$ and $\bbeta$ with $|\alpha_m|, |\beta_n| \leq 1$. Then for any set $\cM\subseteq\mathbb{F}_p^{\times}$ of cardinality $M$ and an interval 
$\cN =   \{B+1, \ldots, B+N\}\subseteq \F_p^\times$ of length $N$ with
$$
M,N\ge p^{\varepsilon},   \quad M^5N^8 \ge p^{5 + \varepsilon},
 \quad MN \ge p^{3/4 + \varepsilon} , \quad M^2N^6 \ge p^{3 + \varepsilon} , 
$$ 
we have
$$
\left|\SABMNp\right| \ll MN p^{-\eta}
$$
where $\eta> 0$ depends only on $\varepsilon$. 
\end{cor} 

 \section{Comparison} 
 First, we note that  it appears that the exponent $1/2 -3/4\ell$ in the 
condition 
$$
N \le \frac{1}{2} q^{1/2 - 3/4\ell}
$$  
of~\cite[Theorem~4.1]{KMS2} can be replaced with  $1/2 +3/4\ell$. 
Now, to see that  Corollary~\ref{cor:SABMN-range} improves on the range 
of~\cite[Theorem~4.1]{KMS2} we choose $N = \fl{p^{3/4 +\delta}}$ 
for some small $\delta>0$. As we have mentioned, in order to apply~\cite[Theorem~4.1]{KMS2} 
one needs to choose $\ell$ with 
$$
N \le \frac{1}{2} q^{1/2 + 3/4\ell}
$$
which leaves only one option $\ell = 2$. 
It is now easy to see that~\cite[Theorem~4.1]{KMS2} is nontrivial only under the condition
$$
\frac{p^{3/4 + 3/8}}{MN} < 1,
$$
which means that it only applies to sets of cardinality $M >  p^{3/8-\delta}$.
On the other hand, for the same choice $N = \fl{p^{3/4 +\delta}}$, 
Corollary~\ref{cor:SABMN-range} 
requires only $M \ge p^{\varepsilon}$ for fixed $\varepsilon >0$.

 \section{Possible Applications} 
As an example of an application of Corollary~\ref{cor:SABMN-range}
we note the bound 
$$
\sum_{m=1}^M\sum_{n=1}^N \alpha_m\beta_n K_{r,p}\(f(m)n\)  \ll MN p^{-\eta}, 
$$
with a non-constant  polynomial $f\in \F_p[X]$
which holds under the same conditions   as in Corollary~\ref{cor:SABMN-range}
and which is unaccessible via the results of~\cite{KMS1,KMS2} unless $M$ is small
and $N$ a rather narrow range. Namely for the bound of~\cite[Theorem~4.1]{KMS2} to apply to the above sum, one needs $f$ is a fixed polynomial, defined over $\Z[X]$ and with bounded coefficients and then one also  needs 
the following rather stringent conditions 
$$
(p^{1+ \varepsilon}/M)^{1/2}  \le N \le p^{1+ \varepsilon}/M^d.
$$
where $d = \deg f$.  Bilinear sums of this type with   quadratic polynomials of the form
$f(X) = (X+A)^2$ (for which our result is still new) appears in applications to the distribution
of squarefree numbers in arithmetic progressions~\cite{Nun}, see 
also~\cite[Section~1.5.1]{KMS1}.

 \section{Mixed additive energy of intervals and arbitrary sets} 
   Given a prime $p$, an integer $H\in [1,p)$, and an arbitrary set $\mathcal{M}\subseteq\mathbb{F}_p^{\times}$,  where $\mathbb{F}_p$, is the finite field with $p$ elements, let $J(H,\mathcal{M})$ denote the number of solution of the congruence 
$$
   xm\equiv ym \pmod p,
$$
   for $x,y\in [1,H)$ and $m,n\in \mathcal{M}$.
   We need the following result of Banks and Shparlinski~\cite[Theorem~2.1]{BaSh}
   
   \begin{lemma}\label{lem:cong}
   For an integer $H\ge 1$ and $\mathcal{M}\subseteq\mathbb{F}_p^{\times}$, of cardinality $M$, the following holds,
\begin{align*}
 J(H,M) & \ll H^2M^2p^{-1}\\ 
&\qquad  + \begin{cases}
HMp^{o(1)}, ~&\text{if}~H\geq p^{2/3};\\
HM^{7/4}p^{-1/4+o(1)}+M^2,~&\text{if}~H< p^{2/3},~M\geq p^{1/3};\\
HMp^{o(1)}+M^2,~&\text{if}~H< p^{2/3},~M< p^{1/3}.
\end{cases}
\end{align*}   
   \end{lemma}
   
We note that in the last case $H< p^{2/3}$ and $M< p^{1/3}$ of Lemma~\ref{lem:cong} the leading term  $H^2M^2p^{-1}$
 never dominates and can be omitted. 
   
\section{Proof of Theorem~\ref{thm:SABMN}} 
   First applying the Cauchy inequality we take the norm of $\{\beta_m\}$ outside
   and obtain 
   \begin{equation}\label{type-II}
   \left|\SABMNp\right|^2\leq \|\beta\|_2^{2}(\|\alpha\|_2^2N+\fS) \leq \|\beta\|_2^{2}(MN+\fS),
   \end{equation}
   where 
$$
\fS=\sum_{\substack{m_1,m_2 \in \cM\\m_1\neq m_2}}\alpha_{m_1}\overline{\alpha_{m_2}}\sum_{n\in \cN}\cK_{r,p}(m_1n)\overline{\cK_{r,p}(m_2n)}.
$$

 We now consider all variables to be elements of $\Fp$ and so we use the language of equations instead
 of congruences modulo $p$. 
 
   Consider two integer parameters $A$ and $B$ such that $2AB\leq N$. Then we follow a strategy followed in \cite[Section 4.1]{KMS2}. In particular,  for some complex weight $\boldsymbol{\eta}=\(\eta_b\)_{b\sim B}$ with $|\eta_b|=1$, to have
\begin{equation}\label{eq-SABMNP}
\fS \ll \frac{\log p}{AB}\sum_{s,t,u\in \mathbb{F}_p^{\times}}\nu(s,t,u)\left|\sum_{b\sim B}\eta_{b}\cK_{r,p}(s(u+b))\overline{\cK_{r,p}(t(u+b))}\right|, 
\end{equation} 
where 
$$
\nu(s,t,u)=\ssssum_{\substack{m_1,m_2 \in \cM,\ m_1\neq m_2, \\ n \in \cN\
a\sim A\\ am_1=s,\ am_2=t, \ \overline{a}n=u}}\left|\alpha_{m_1}\alpha_{m_2}\right|.
$$
Now,  writing 
$$
\nu(s,t,u) = \nu(s,t,u)^{(\ell-1)/\ell}   (\nu(s,t,u)^2)^{1/2 \ell}
$$
and then using the H\"{o}lder inequality in~\eqref{eq-SABMNP} with weights $\ell/(\ell-1)$, $2\ell$ and $2\ell$, we have
\begin{equation}\label{G}
\left|\fS \right|^{2\ell}\leq \frac{1}{(AB)^{2\ell}} R_1^{2\ell-2}R_2\SKr p^{o(1)},
\end{equation}
where 
$$
R_1=\sum_{s,t,u\in \mathbb{F}_p^{\times}}\nu(s,t,u)
\mand
R_2=\sum_{s,t,u\in \mathbb{F}_p^{\times}}\nu(s,t,u)^2,
$$
and
$$
\SKr=\sum_{\substack{s,t,u\in \F_p^{\times}\\s\ne t}}\left|\sum_{b\sim B}\eta_{b}\cK_{r,p}(s(u+b))\overline{\cK_{r,p}(t(u+b))}\right|^{2\ell}.
$$
Trivially, we have
$$
\sum_{s,t,u\in \mathbb{F}_p^{\times}}\nu(s,t,u)\ll AM^2N.
$$
Consider the set $\mathcal{A}$, which contains integers $a$ such that $a\sim A$. Now using the fact that $|\alpha_m|\leq 1$, we have
 \begin{align*}
 R_2 &\leq \Bigl\{\(a_1,a_2,m_1,m_2,m_3,m_4,n_1, n_2\)\in\mathcal{A}^2\times \mathcal{M}^4\times\mathcal{N}^2:\\
 & \qquad \qquad \qquad \qquad \overline{a_1}n_1=\overline{a_2}n_2, \ a_1m_1=a_2m_2, \ a_1m_3=a_2m_4\Bigr\}.
 \end{align*}
 Now if we fix  one of $J(2A,\mathcal{M})$  solutions, $(a_1,a_2,m_1,m_2)\in \mathcal{A}^2\times \mathcal{M}^2$ to $a_1m_1=a_2m_2$, then there are at most $MN$ many solutions
  $(m_3, m_4, n_1,n_2)\in \cM^2\times \cN^2$  to
 $$
 \overline{a_1}n_1=\overline{a_2}n_2 \mand a_1m_3=a_2m_4.
$$
 Hence the above bound becomes,
$$
 R_2\ll MN J(2A,\mathcal{M}).
$$
 Now to bound $\SKr $, we use a result of Kowalski, Michel and Sawin~\cite[Equation~(4.7)]{KMS2}, which together with the bounds
 of~\cite[Theorem~4.5]{KMS2} gives 
$$
   \SKr \ll p^3B^{\ell}+p^2B^{2\ell- \rf{(\ell-1)/2}}+p^{3/2}B^{2\ell}.
$$
  Choose $B=\fl{0.25 p^{3/2\ell}}$.  
  This gives
$$
  \SKr \ll p^3B^{\ell}+p^2B^{2\ell- \rf{(\ell-1)/2}}.
$$
We now observe that for $\ell \ge 2$ we have 
$$
\rf{(\ell-1)/2} \ge \ell/3. 
$$
Hence for the above choice of $B$ we have 
$$
p^2B^{2\ell- \rf{(\ell-1)/2}} \le p^2B^{5\ell/3}=   p^{9/2} =  p^3B^{\ell}.
$$
Therefore, 
$$
  \SKr \ll p^3B^{\ell}.
$$
and  from~\eqref{G}, we have
\begin{align*}
\left|\fS \right|^{2\ell}&\leq \frac{1}{(AB)^{2\ell}} (AM^2N )^{2\ell-2}J(2A,\mathcal{M})MN B^{\ell}p^{3+o(1)}\\
&=A^{-2}B^{-\ell}M^{4\ell-3}N^{2\ell-1}J(2A,\mathcal{M})p^{3+o(1)}.
\end{align*}
Now we apply Lemma~\ref{lem:cong}  for the bound for $J(2A,M)$ and get
\begin{align*}
\left|\fS \right|\leq A^{-1/\ell}&M^{2-3/{2\ell}}N^{1-1/{2\ell}}\\
&\(A^2M^2p^{-1}+AM^{7/4}p^{-1/4}+AM+M^2\)^{1/{2\ell}}p^{{3/{4\ell}}+o(1)}.
\end{align*} 
Choosing $A =  \fl{0.5 NB^{-1}} $, 
so that 
$$
N p^{-3/2\ell} \ll A \ll N p^{-3/2\ell},
$$
 we  now derive 
 $$
 \left|\fS \right|\leq M^2N \Gamma^{1/2\ell} p^{o(1)}, 
 $$
 where 
\begin{align*}
\Gamma& \leq A^{-2}M^{-3}N^{-1}\(A^2M^2p^{-1}+AM^{7/4}p^{-1/4}+AM+M^2\)p^{3/{2}}\\
&=M^{-1}N^{-1}p^{1/2}+A^{-1}  M^{-5/4}N^{-1} p^{5/4}\\
& \qquad  \qquad +A^{-1} M^{-2} N^{-1}  p^{3/2} +A^{-2} M^{-1} N^{-1} p^{3/2} \\
&=M^{-1}N^{-1}p^{1/2}+  M^{-5/4} N^{-2} p^{5/4+3/2\ell}\\
& \qquad  \qquad +  M^{-2} N^{-2} p^{3/2+3/2\ell} + M^{-1} N^{-3} p^{3/2+3/\ell}.
\end{align*}
Inserting this bound in~\eqref{type-II}, after simple calculations, we complete the proof.

\section{Comments} 

It is easy to see that under the Generalised Riemann Hypothesis (GRH), the bound of  Lemma~\ref{lem:cong}
can be improves as 
$$
 J(H,M) =  H^2M^2p^{-1} + O\(HMp^{o(1)}\), 
$$
see~\cite[Equation~(1.4)]{BaSh}. In turn,  this implies that under the GRH, one can recover all 
main results of~\cite{KMS2} for arbitrary sets without any restrictions on 
$$M^{+} = \max\{m:~m \in \cM\}, $$
which are imposed  in~\cite{KMS2}.  

We also use the opportunity to note that using the results of~\cite[Theorem~1]{ACZ},
one can extend~\cite[Theorems~4.1 and~4.3]{KMS2} to intervals $\cN$ in  an arbritry 
position rather than just at the origin. 

On the other hand, the method of~\cite{KMS1} 
uses completing technique and thus it is not clear how to extend it to arbitrary sets $\cM$. 

\section*{Acknowledgement}

During the preparation of this work
N.B. was supported by the post-doctoral fellowship, Harish-Chandra Research Institute, Prayagraj, India,  and 
 I.E.S.   by the  Australian Research Council  Grant~DP170100786
and by the Natural Science Foundation of China Grant~11871317.

\end{document}